\documentclass[aip,cha,epsf,reprint,groupedaddress]{revtex4-1}
\usepackage{amsmath,amsthm,amssymb,amstext,amsfonts}
\usepackage{graphicx}
\usepackage{bm}
\usepackage{comment}
\usepackage{placeins}
\usepackage{array}
\usepackage{caption}
\usepackage{subcaption}
\usepackage{color} 
\usepackage{hyperref}
\usepackage{soul,xcolor}
\setstcolor{blue}

\newcommand{\be}{\begin{enumerate}}
\newcommand{\ee}{\end{enumerate}}
\newcommand{\mrm}{\mathrm}

\newcommand{\etal}{{\em et al.}~}

\newcommand{\R}{\mathbb{R}}

\newcommand{\E}{\mrm{E}}

\begin{document}
\title{Geometric k-nearest neighbor estimation of entropy and mutual information}
\author{Warren M. Lord}
\author{Jie Sun}
\author{Erik M. Bollt}
\affiliation{Department of Mathematics, Clarkson University, 8 Clarkson Ave, Potsdam, NY, 13699-5815, USA}
\date{\today}

\begin{abstract}
\textcolor{black}{Nonparametric estimation of mutual information is used in a wide range of scientific problems to \textcolor{black}{quantify} dependence between variables. The $k$-nearest neighbor (knn) methods are consistent, and therefore expected to work well for large sample size. These methods use geometrically regular local volume elements. This practice allows maximum localization of the volume elements, but can also induce a bias due to a poor description of the local geometry of the underlying probability measure.} We introduce a new class of knn estimators that we call geometric knn estimators (g-k\textcolor{black}{nn}), which use more complex local volume elements to better model the local geometry of the probability measures. As an example of this class of estimators, we develop a g-k\textcolor{black}{nn} estimator of entropy and mutual information based on elliptical volume elements, capturing the local stretching and compression common to a wide range of dynamical systems attractors. \textcolor{black}{A} series of numerical examples \textcolor{black}{in which the thickness of the underlying distribution and the sample sizes are varied suggest that local geometry is a source of problems for knn methods such as} the Kraskov-St\"{o}gbauer-Grassberger (KSG) estimator \textcolor{black}{when local geometric effects cannot be removed by global preprocessing of the data.} \textcolor{black}{The g-knn method performs well despite the manipulation of the local geometry.} In \textcolor{black}{addition}, the examples suggest that the g-k\textcolor{black}{nn} estimators can be of particular relevance to applications in which the system is large\textcolor{black}{,} but data size is limited.
\end{abstract}


\maketitle

\begin{quotation}
  Mutual information is a tool used by scientists to quantify the dependence between variables without making specific modeling assumptions. \textcolor{black}{In many applications mutual information must be estimated from a finite set of data with no model specific assumptions about the distributions of the variables.} The most popular estimators of mutual information are $k$-nearest neighbor (knn) estimators, which \textcolor{black}{locally estimate the distributions from statistics of distances between data points.} \textcolor{black}{The $k$nn methods typically} use \textcolor{black}{geometrically regular} local volume elements. We introduce a new class of knn estimators, the geometric knn estimators (g-k\textcolor{black}{nn}), that use more detailed local volume elements to model the geometric features of the probability density function. Inspired by the local geometry of dynamical systems attractors, we develop a singular value decomposition driven g-k\textcolor{black}{nn} estimator that models local volumes by ellipsoids. We show by numerical examples that g-k\textcolor{black}{nn} estimators can alleviate bias due to \textcolor{black}{thinly supported distributions and} small sample size.
\end{quotation}

\section{Introduction}

Differential mutual information is used in a number of scientific disciplines to measure the strength of the relationship between two continuous random variables. Given two continuous random variables, $X$ and $Y$, taking values in $\R^{d_X}$ and $\R^{d_Y}$, the differential entropy, $H$, and mutual information, $I$ are defined by
 \begin{align}
   H(X)&=-\E[\log(f_X(X))]  \label{eq:defentropy}\\
   I(X;Y) &= H(X)+H(Y)-H(X,Y),\label{eq:defmutualinformation}
 \end{align}
 where $f_X$ is the probability density function (pdf) of $X$, $\E[\cdot]$ is the expected value functional, and $H(X,Y)$ is the entropy of the joint variable $(X,Y)$. \textcolor{black}{In dynamical systems applications the variables $X$ and $Y$ are often produced by high dimensional nonlinear dynamical or stochastic process so that exact computation of Eqs.~\eqref{eq:defentropy} and \eqref{eq:defmutualinformation} is impractical. Furthermore, when data is generated empirically, the model is often unknown, and there are many applications\cite{sun2015causal,lord2016inference} in which a large number of such estimates must be made in an automated manner. Therefore, the problem of nonparametric estimation of differential entropy and mutual information, in which $N$ points in $\R^d$ are used to estimate Eqs.~\eqref{eq:defentropy} or \eqref{eq:defmutualinformation} without model specific assumptions, has received much attention from the statistics, probability, machine-learning, and dynamics communities~\cite{stowell2009fast,darbellay1999estimation,beirlant1997nonparametric,joe1989estimation,kozachenko1987sample,KSG,calsaverini2009information,giraudo2013non}.}

\begin{figure}[]
  \centering
  \begin{subfigure}[b]{.57\columnwidth}
    \includegraphics[width=.9\linewidth]{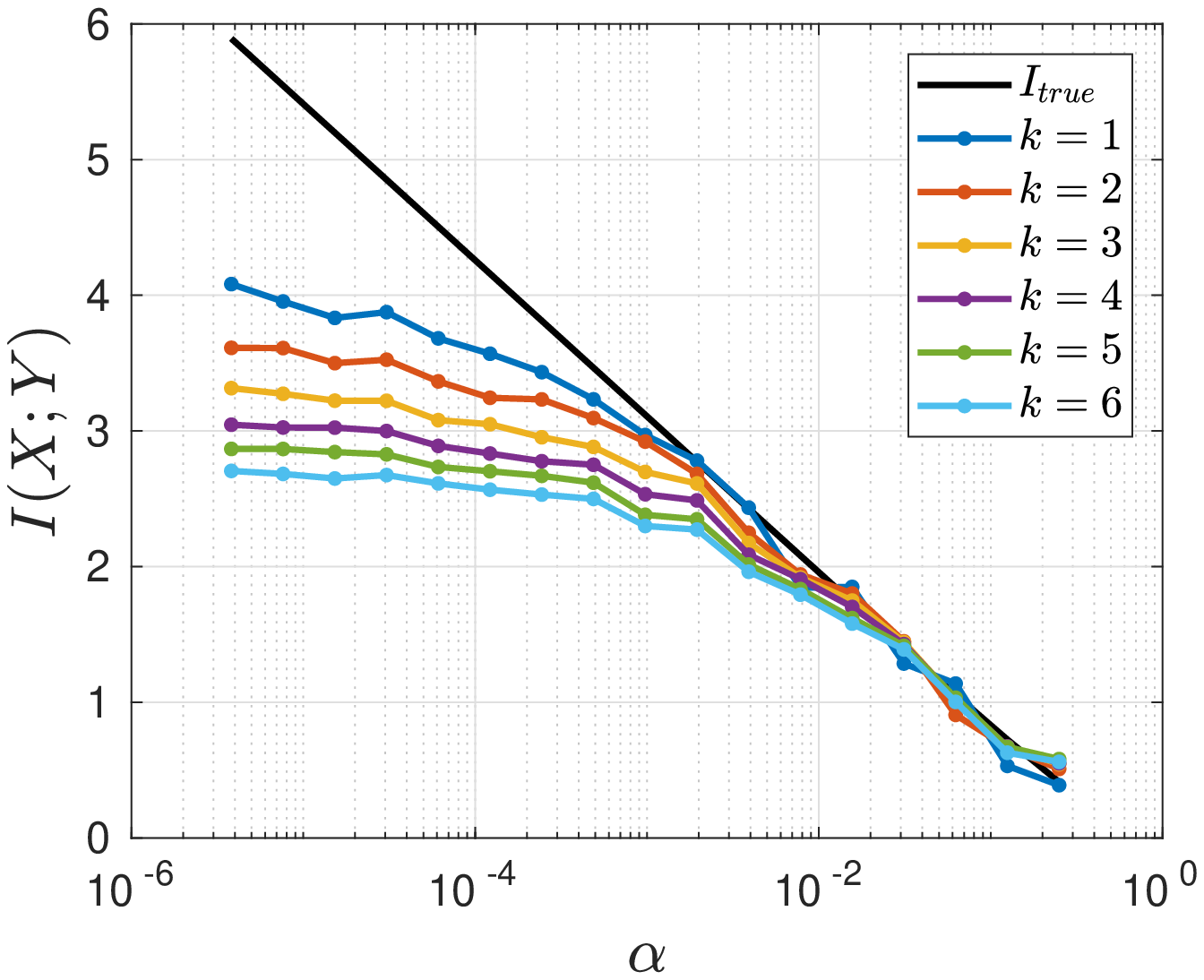}
    \caption{}
    \label{fig:ksgfail}
  \end{subfigure}\hfill
  \begin{subfigure}[b]{.42\columnwidth}
    \centering
    \includegraphics[width=\linewidth]{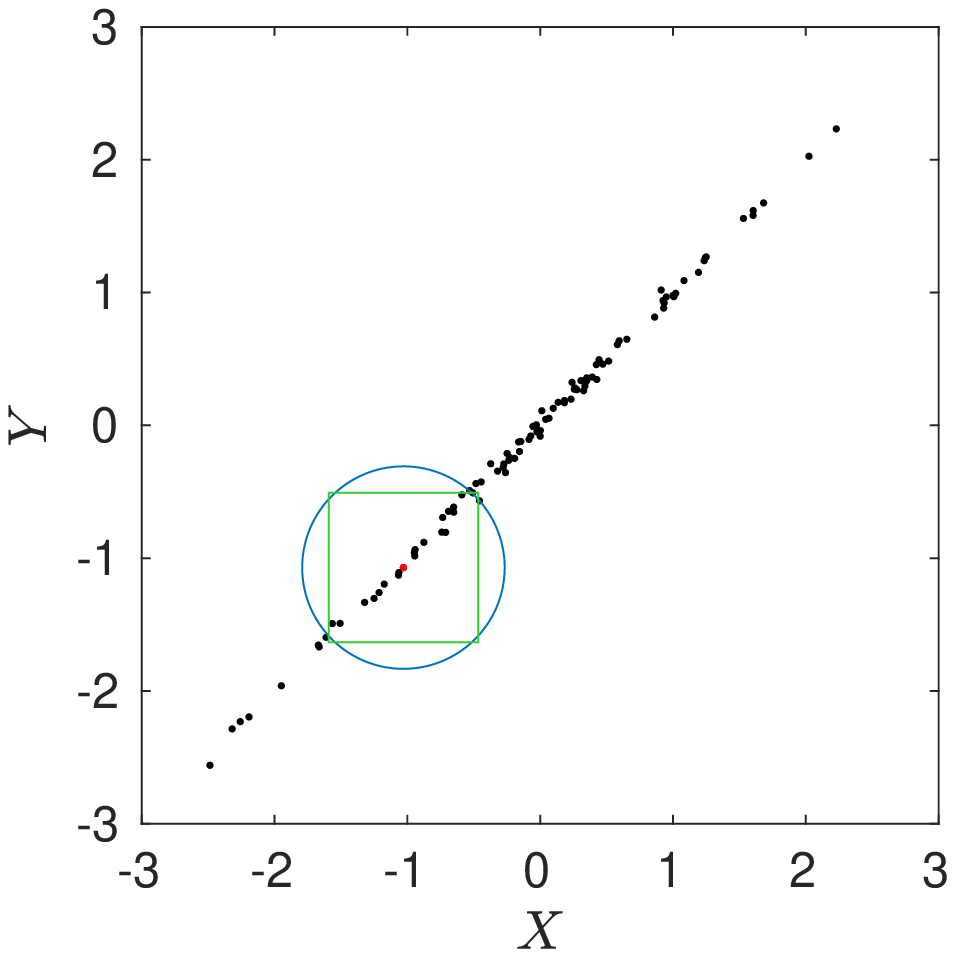}
    \caption{}
    \label{fig:sphericallinear}
  \end{subfigure}
    \caption{\textcolor{black}{(a) KSG estimates of mutual information for two 1d normally distributed variables with standard deviation 1 and correlation $1-\alpha$. For each $\alpha\in\{2^{-j}:j=2,\ldots,18\}$ a sample of size $N=100$ is drawn and the mutual information estimated by KSG with $k=1,\ldots,6$. The true mutual information, $I_{true}=I(X;Y)$ is plotted in black. (b) A sample of size 100 when $\alpha=.001$. A randomly chosen sample point is highlighted in red. A sphere in the maximum norm is plotted in green and a sphere in the Euclidean norm is plotted in blue. The radius of each sphere is equal to the distance to the 20th closest neighbor in the respective norm.}}
  \label{fig:mvnormalexample}
\end{figure}

 The \textcolor{black}{$k$-nearest neighbors (}$k$nn\textcolor{black}{)} estimators have received particular attention due to their ease of implementation and efficiency in a multidimensional setting. \textcolor{black}{The $k$nn estimates
 are derived by expressing volumes of neighborhoods of data points in terms of distances from each data point to other data points in the neighborhood. In most $k$nn methods these local volume elements are assumed to be highly regular in order to minimize the amount of data needed to define them, and therefore keep the volume elements as local as possible. For instance, the Kozachenko\textcolor{black}{-}Leonenko estimator for differential entropy~\cite{kozachenko1987sample} uses $p$-spheres as local volume elements. The popular Kraskov-St\"{o}gbauer-Grassberger (KSG) estimator for differential mutual information is based on \textcolor{black}{Kozachenko-Leonenko} estimates of each of the entropies in~\eqref{eq:defmutualinformation} in which the local volume elements in the estimate of $H(X,Y)$ are taken to be products of the volume elements in the marginal distributions in order to achieve cancellations of the bias. \textcolor{black}{The estimator has} two versions, one in which all volume elements are max-norm spheres, such as illustrated by the green max-norm sphere in Fig.~\ref{fig:sphericallinear}, and the other in which the volume elements in the joint space are the product of $p$-norm spheres in the marginal spaces.}

\textcolor{black}{The most important advantage to using such $k$nn methods based on highly regular volume elements is their asymptotic consistency~\cite{singh2003nearest,gao2017demystifying}. A drawback in} data-driven applications \textcolor{black}{where} sample size is fixed and often limited is that the local volume elements might not be descriptive of the geometry of the underlying probability measures, resulting in bias in the estimators. A simple example of this problem is shown in Figure~\ref{fig:ksgfail}, in which $X$ and $Y$ are normally distributed with standard deviation $1$ and correlation $1-\alpha$. By direct computation\textcolor{black}{,} the true mutual information increases asymptotically as $\log(\alpha)$, but for each $k$ the KSG estimator applied to the raw data diverges quickly as $\alpha$ decreases. Figure~\ref{fig:sphericallinear} describes the idea that the problem may be due to local volume elements not being descriptive of the geometry of the underlying measure. Improving on that issue is the major stepping off point of this paper. In particular, the KSG local volume elements mostly resemble the green square (a max-norm sphere), whose volume greatly overestimates the volume spanned by the data points it contains.

\textcolor{black}{This paper introduces a new class of $k$nn estimators, the g-knn estimators, which use} more irregular \textcolor{black}{local volume elements} that are more descriptive of the underlying geometry at the smallest length scales represented in the data. \textcolor{black}{The defining feature of g-knn methods} is a trade-off between the irregularity of the object, which requires more local data to fit\textcolor{black}{, and therefore less localization of the volume elements}, and the improvement in the approximation of the \textcolor{black}{local geometry of the} underlying measure. Because of this trade-off, the local volume element should be chosen to reflect the geometric properties expected in the desired application. 

Motivated by the study of dynamical systems, it is reasonable to model these properties on the local geometry of attractors. One of the most striking geometric features of attractors common to a wide class of dynamical systems, including dissipative systems and dynamical systems with competing time scales, is stretching and compression in transverse directions~\cite{ott2002chaos}. More precisely, the local geometry is characterized by both positive Lyapunov exponents corresponding to directions in which nearby points are separated over time and negative Lyapunov exponents corresponding to orthogonal directions in which nearby points are compressed.

\textcolor{black}{To test the idea behind the g-knn estimators, Sec.~\ref{sec:method} develops a particular g-knn method that uses local volume elements to match the geometry of stretching and compression in transverse directions.} Ellipsoids are a good option for capturing this geometric feature because they have a number of orthogonal axes with different lengths. They are also fairly regular geometric objects: the only parameters that require fitting are the center and one \textcolor{black}{axis} for each dimension. Such ellipsoids can be fit very efficiently using the singular value decomposition (svd) of a matrix formed from the local data~\cite{golub2012matrix}. 

\textcolor{black}{The g-knn estimator is tested on four} one-parameter families of joint random variables \textcolor{black}{in which the parameter controls the stretching of the geometry of the underlying measure}. \textcolor{black}{The estimates} are \textcolor{black}{compared} with the KSG estimator as the local geometry of the joint distribution becomes more stretched. \textcolor{black}{D}istributions can \textcolor{black}{also} appear to be more stretched locally \textcolor{black}{if} local neighborhoods of data increase in size, which occurs in $k$nn methods when sample size is decreased. Therefore, the g-k\textcolor{black}{nn} estimator and KSG are also compared numerically in examples using small sample size.

Unlike the \textcolor{black}{Kozachenko-Leonenko} and KSG estimators, the g-knn method developed here has not been corrected for asymptotic bias, so that it should be expected that KSG outperforms this particular g-knn method for large sample size. \textcolor{black}{What is surprising is that the g-knn estimator developed in Sec.~\ref{sec:method} outperforms KSG for small sample sizes and thinly supported distributions despite lacking the clever bias cancellation scheme that defines KSG.}  Since KSG is considered to be state-of-the-art in the nonparametric estimation of mutual information, the result should hold for other methods that do not account for local geometric effects.

There have been many attempts to resolve the bias of KSG. For instance, Zhu \etal~\cite{zhu2014bias} improved on the bias of KSG by expanding the error in the estimate of the expected amount of data that lies in a local volume element. Also, Wozniak and Kruszewski~\cite{wozniak2012estimating} improved KSG by modeling deviations from local uniformity using the distribution of local volumes as $k$ is varied. The\textcolor{black}{se} improvements do not \textcolor{black}{directly} address the limitations of spheres to describe interesting features of the local geometry.

The \textcolor{black}{class} of g-k\textcolor{black}{nn} estimator\textcolor{black}{s} can be thought of as generalizing the estimator of mutual information described by Gao, Steeg, and Galstyan (GSG) in 2015\cite{gao2015efficient}, which uses a principle component analysis of the local data to fit a hyper-rectangle. The \textcolor{black}{svd-based} g-k\textcolor{black}{nn} estimator \textcolor{black}{defined in Sec.~\ref{sec:method} improves on the} GSG \textcolor{black}{treatment of local data}. \textcolor{black}{These improvements are highlighted in} Sec.~\ref{sec:method}.

\section{Method}\label{sec:method}

This section defines a g-k\textcolor{black}{nn} estimate of entropy that\textcolor{black}{, in turn,} yields an estimate of mutual information when substituted in Eq.~\eqref{eq:defmutualinformation}. Let the given data set be denoted $\{x_i\subset \R^d:i=1,\ldots,N\}$, where for each $i$, $x_i$ is a sample point. The g-k\textcolor{black}{nn} estimate of entropy is similar to \textcolor{black}{the} Kozachenko\textcolor{black}{-}Leonenko estimator for entropy \cite{kozachenko1987sample} in that the entropy is estimated using the resubstitution formula
\begin{align}
  H(X)&=-\E[\log f_X(X)] \notag\\
   &\approx -\frac1N\sum_{i=1}^N\log(\widehat{f}_X(x_i)),\label{eq:resub}
\end{align}
where $\widehat{f}_X(x_i)$ is an estimate of the pdf $f_X$ at $x_i$.
The pdf $f_X$ at $x_i$ is estimated by
\begin{align}
  \widehat{f}_X(x_i)\propto \frac{k(x_i)/N}{\mrm{Vol}_i},\label{eq:volfrac}
\end{align}
where $k(x_i)$ is the number of data points other than $x_i$ in a neighborhood of $x_i$, and $\mrm{Vol}_i$ is the volume of the neighborhood.

\subsection{SVD estimation of volume elements}
The elliptical local volume elements are estimated by singular value decomposition (svd) of the local data.
For any fixed $x_i$, denote the $k$ nearest neighbors by the vectors $x_i^j$, $j=1,\ldots,k$ in $\R^d$. Define the $k$-neighborhood of sample points to be $\{x_i^j:j=0,\ldots,k\}$ where $x_i^0=x_i$. In order for the svd to indicate directions of maximal stretching it is first necessary to center the data. Let $z=\frac{1}{k+1}\sum_{j=0}^k x_i^j$ be the centroid of the neighborhood in $\R^d$, and define the centered data vectors in $\R^d$ by $y^j=x_i^j-z$. In order for the svd, a matrix decomposition, to operate on the centered data, form the $(k+1)\times d$ matrix of row vectors,
\begin{align}
  Y=\left(\begin{array}{c} y^0 \\ y^2 \\\vdots \\ y^{k}\end{array}\right).
\end{align}

Since $Y$ is a $(k+1)\times d$ matrix it has an svd of the form $Y=U\Sigma V^T$, where $U$ is a $(k+1)\times(k+1)$ unitary matrix, $\Sigma$ is a $(k+1)\times d$ dimensional matrix which is zero with the possible exception of the nonnegative diagonal components, and $V$ is a $d\times d$ unitary matrix. The columns of $U$ and $V$ are the left and right singular vectors, and since the left singular vectors do not play a role in this estimator, the word ``right'' will be omitted when referring to the right singular vectors.

Since $V$ is unitary, the singular vectors, $v_i^{(l)}$ are of unit length and orthogonal. The first singular vector, $v_i^{(1)}$ points in the direction in which the data is stretched the most, and each subsequent singular vector points in the direction \textcolor{black}{which is} orthogonal to all previous singular vectors \textcolor{black}{and} that accounts for the most stretching.

The singular values are equal to the square root of the sum of squares of the lengths of the projections of the $y^j$ onto a singular vector. This means that $\sigma_i^l$ is $\sqrt{k}$ times the standard deviation of the projection of the data onto $v_i^{(l)}$.

The use of data centered to the mean is an important difference between the g-k\textcolor{black}{nn} estimator described here and the \textcolor{black}{GSG} estimator, in which the data is centered to $x_i$ \textcolor{black}{(see Footnote~2 in Ref.~\cite{gao2015efficient})}. Centering the data at $x_i$ can bias the direction of the singular vectors away from the directions implied by the underlying geometry. In Fig.~\ref{fig:neighborhoodellipse}, for instance, the underlying distribution from which the data is sampled can be described as constant along lines parallel to the diagonal $y=x$ and a bell curve in the orthogonal direction, with a single ridge along the line $y=x$. If local data were centered at the red data point then all vectors would have positive inner product with \textcolor{black}{the vector} $(-1,1)$, so that the first singular vector would be biased toward $(-1,1)$. The center of the local data, on the other hand, is near the top of the ridge, so that the singular vectors of the centered data (in blue in the figure) estimate the directions along and transverse to this ridge.

\subsection{Translation and scaling of volume elements}
\begin{figure}
  \includegraphics[width=.49\textwidth]{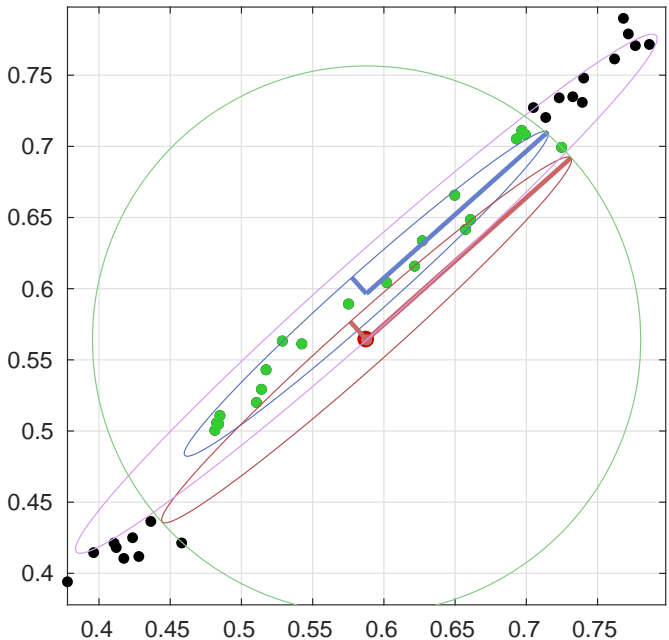}
  \caption{A sample from a random variable with one sample point, $x_i$, highlighted in red at approximately $(0.59,0.56)$, and its $k=20$ nearest neighbors in Euclidean distance highlighted in green. The ellipsoid centered at $x_i$ and drawn in red contains the volume used by the g-K\textcolor{black}{nn} estimator. The length of the major axes is determined by the largest projection of one of the $k$ neighbors onto the major axes, enclosing $k_i=3$ points in its $k$-\textcolor{black}{neighborhood}, including itself. Two other ellipsoids are centered at the centroid of the $k+1$ neighbors. The larger one (magenta) has radii large enough to enclose all $k+1$ data points. The major axis of the smaller ellipsoid (blue) is determined by the largest projection of a data point onto the major axis. All three ellipsoids have the same ratio of lengths of axes, $\sigma_1/\sigma_2$, where the $\sigma_i$ are the singular values determined by the centered $k$-\textcolor{black}{neighborhood}.}
  \label{fig:neighborhoodellipse}
\end{figure}

Since the $v_i^{(l)}$ are orthogonal, the vectors $\sigma_i^lv_i^{(l)}$ can be thought of as the axes of an ellipsoid centered at the origin. The ellipsoid needs to be translated to the $k$-neighborhood and scaled to fit the data. There are many ways to perform this translation and rescaling, three of which are depicted in Fig.~\ref{fig:neighborhoodellipse}.

In Fig.~\ref{fig:neighborhoodellipse}, there are two ellipsoids centered on the centroid of the $k$-neighborhood. The larger ellipsoid is the smallest ellipsoid that contains or intersects all points in the $k$-neighborhood. The problem with this approach is that these ellipsoids might contain data points which are not one of the $k$ nearest neighbors of $x_i$, as is seen in Fig.~\ref{fig:neighborhoodellipse}. One solution to this problem is to include these data points in the calculation of the proportion of the data that lies inside the volume defined by the ellipsoid. We avoid this approach, however, both because it involves extra computational expense in finding these points, and because the new neighborhood obtained by adjoining these points is less localized.

An alternative approach is to decrease the size of the ellipsoid to exclude points not in the $k$-neighborhood, as depicted by the smaller ellipsoid centered at the centroid. Such an ellipsoid could contain a proper subset of the $k$-neighborhood so that $k(x_i)$ may be less than $k$. The problem with this approach, however, is that in higher dimensions, the ellipsoid may contain no data points, introducing a $\log(0)$ into Equation~\eqref{eq:resub}.

Instead of centering at the centroid, however, the ellipsoid could be centered at $x_i$. If the length of the major axis is taken to be the Euclidean distance to the furthest neighbor in the $k$-neighborhood, then the ellipsoid and its interior will only contain data points in the $k$-neighborhood because the distance from $x_i$ to any point on the ellipsoid is less than or equal to the distance to the furthest of the $k$ neighbors \textcolor{black}{(the sphere in Fig.~\ref{fig:neighborhoodellipse})}, which is less than or equal to the distance to any point not in the $k$-neighborhood. This neighborhood will contain at least one data point. An example of such an ellipsoid is shown in Fig.~\ref{fig:neighborhoodellipse}, which includes $k(x_i)=3$ data points.

The \textcolor{black}{GSG} estimator\cite{gao2015efficient} is centered at $x_i$, but the lengths of the sides of the hyper-rectangles seem to be determined by the largest projection of the local data onto \textcolor{black}{the axes}, which destroys the ratio of singular vectors that describes the local geometry. In addition, it is possible that the corners of the hyper-rectangles may circumscribe more than the $k$ neighbors of $x_i$, even though a constant $k$ neighbors are assumed in the estimate.

\subsection{The g-k\textcolor{black}{nn} estimates for $H(X)$ and $I(X;Y)$}
In order to explicitly define the global estimate of entropy that results from this choice of center, define $\epsilon(x_i,k)$ to be the Euclidean distance from $x_i$ to the $k$th closest data point. \textcolor{black}{Define}
\begin{align}
  r_i^l &=\epsilon(x_i,k)\frac{\sigma_i^l}{\sigma_i^1}
\end{align}
\textcolor{black}{to be the lengths of the axes of the ellipsoid centered at $x_i$, for $l=1,\ldots,d$.}
Note that $r_i^1=\epsilon(x_i,k)$, and
\begin{align}
  \frac{r_i^{l_1}}{r_i^{l_2}}&=\frac{\sigma_i^{l_1}}{\sigma_i^{l_2}}.
\end{align}

The volume of this ellipsoid can be determined from the formula
\begin{align}
  V_i &= \frac{\pi^{d/2}}{\Gamma(1+\frac{d}{2})}\prod_{l=1}^d r^l_i \\
      &= \frac{\pi^{d/2}}{\Gamma(1+\frac{d}{2})}\epsilon(x_i,k)^d\prod_{l=1}^d \frac{\sigma^l_i}{\sigma_i^1}.
\end{align}

Substitution into Equations~\eqref{eq:resub} and~\eqref{eq:volfrac} yields
\begin{widetext}
\begin{align}
  \widehat{H}_{g-k\textcolor{black}{nn}}(X) =& -\frac1N \sum_{i=1}^N \log\frac{k(x_i)\Gamma(1+\frac{d}{2})}{N\pi^{d/2}\epsilon(x_i,k)^d\prod_{l=1}^d (\sigma_i^l/\sigma_i^1)} \\
  =& \log(N)+\log(\pi^{d/2}/\Gamma(1+d/2)) - \frac1N\sum_{i=1}^N\log(k(x_i))  
   +\frac{d}{N}\sum_{i=1}^N\log(\epsilon(x_i,k)) +  \frac1N\sum_{i=1}^N\sum_{l=1}^d\log \left(\frac{\sigma_i^l}{\sigma_i^1}\right)
\end{align}
\end{widetext}
The estimate for $I(X;Y)$ is then obtained using Eq.~\eqref{eq:defmutualinformation}. The term $\frac1N\sum_{i=1}^N\sum_{l=1}^d\log \left(\frac{\sigma_i^l}{\sigma_i^1}\right)$ is small when the local geometry is relatively flat, but, as is demonstrated in Section~\ref{sec:examples}, it can have a large impact on the estimate for more interesting local geometries.

\section{Examples}\label{sec:examples}

This section compares KSG estimates of mutual information on simulated examples with the estimates of the g-k\textcolor{black}{nn} estimator defined in Sec~\ref{sec:method}. \textcolor{black}{The examples are designed so that the local stretching of the distribution is controlled by a single scalar parameter, $\alpha$. Plotting estimates against $\alpha$ suggests that the local stretching is a source of bias for KSG, but that the g-knn estimator is not greatly affected by the stretching.}

\textcolor{black}{The} examples \textcolor{black}{are divided into four} one-parameter families of distributions in which the parameter \textcolor{black}{$\alpha$} affects local geometry. \textcolor{black}{Each family is defined by a model, consisting of the distributions of a set of variables and the equations that describe how these variables are combined to create $X$ and $Y$. The objective is to estimate $I(X;Y)$ directly from a sample of size $N$ without any knowledge of the form of the model. }In the first set of examples \textcolor{black}{(Sec.~\ref{sec:closedform})} the \textcolor{black}{models are simple enough that the} mutual information can be computed exactly. In the second set of examples \textcolor{black}{(Sec.~\ref{sec:henon})}, which are designed to be more typical of dynamical systems research, the mutual information of a pair of coupled H\'{e}non maps is estimated for varying amounts of noise. \textcolor{black}{In the latter case t}he qualitative behaviors of the estimators are compared since the system is too complicated to find the true mutual information.

\subsection{Tests where mutual information is known}\label{sec:closedform}
This section defines three one-parameter famil\textcolor{black}{ies} of distributions in which the parameter describes the ``thickness'' of the distribution. \textcolor{black}{Families 1 and 2 are 2d examples with 1d marginals built around the idea of sampling from a 1d manifold with noise in the transverse direction.} The third \textcolor{black}{family} is a 4d joint distribution with 2d marginals.

\begin{description}
\item{\textbf{Family 1 : }}\textcolor{black}{The model is
    \begin{align}
      Y&=X+\alpha V \\
      X,V &\mrm{\ i.i.d.}\, Unif(0,1).
    \end{align}
    }
The result is a support that is a thin parallelogram around the diagonal $Y=X$. \textcolor{black}{As $\alpha\to 0$ the distributions become more concentrated around the diagonal $Y=X$.} The mutual information of $X$ and $Y$ is
  \begin{equation}
    I(X;Y)=-\log(\alpha)+\alpha-\log(2).
  \end{equation}

\item{\textbf{Family 2 : }}The second example is meant to capture the idea that noise usually has some kind of tail behavior. In this example $V$ is a standard normal, so that the noise term, $\alpha V$ is normally distributed with standard deviation $\alpha$. \textcolor{black}{The model is
    \begin{align}
      Y&=X+\alpha V \\
      X&\sim Unif(0,1),\\
      V&\sim\mathcal{N}(0,1),
    \end{align}
    where $X$ and $V$ are independent.
    }
  In this case \textcolor{black}{the exact form of the mutual information is}
  \begin{align}
    I(X;Y)=&-\log(\alpha)+\Phi\left(-\frac{y}{\alpha}\right)\notag\\
     &-\Phi\left(\frac{1-y}{\alpha}\right)-\frac12\log(2\pi e),
  \end{align}
  where $\Phi$ is the cdf of the standard normal distribution.

\item{\textbf{Family 3 : }}In the third example the joint variable is distributed as \textcolor{black}{
    \begin{align}
      (X,Y)&\sim\mathcal{N}(0,\Sigma)\\
    \Sigma&=
  \left[
  \begin{array}{c c | c c}
    7 & -5 & -1 & -3 \\
    -5 & 5 & -1 & 3 \\\hline
    -1 & -1 & 3 & -1 \\
   -3 & 3 & -1 & 2+\alpha
  \end{array}
                  \right].
    \end{align}
    }
The first two coordinates of this variable belong to the variable $X$ and the third and fourth to the variable $Y$. Thus, the upper left 2 by 2 block is the covariance matrix of $X$ and the bottom 2 by 2 block is the covariance of $Y$. As long as $\alpha>0$, $\Sigma$ is positive definite but if $\alpha=0$ then $\Sigma$ is not of full rank, and the distribution $\mathcal{N}(0,\Sigma)$ is called degenerate, and is supported on a 3d hyperplane. When $\alpha$ is positive but small, the distribution can be considered to be concentrated near a 3d hyperplane. \textcolor{black}{In this case t}he mutual information of $X$ and $Y$ is
\begin{equation}
  I(X;Y)=-\frac12\log\left(\frac{|\Sigma_X||\Sigma_Y|}{\textcolor{black}{|}\Sigma\textcolor{black}{|}}\right)
\end{equation}
\end{description}

\textcolor{black}{These examples are simple enough that, instead of using g-knn, one might be able to guess the algebraic form of the model and perform preprocessing to isolate noise and consequently remove much of the bias due to local geometry. For Families 1 and 2, for instance, if the algebraic form of the model $Y=X+\alpha V$ is known, one can express the mutual information as $I(X;Y)=H(Y)-H(\alpha V)+I(X;\alpha V)$, where the term $I(X;\alpha V)$ can be estimated by KSG after dividing by standard deviations (or, if $\alpha V$ is thought to be independent noise, then it would be assumed that $I(X;\alpha V)=0$). This rearrangement of variables is in essence a global version of what the g-knn estimator accomplishes locally using the svd.}

\begin{figure*}[]
  \begin{subfigure}[b]{.37\textwidth}
    \includegraphics[width=\textwidth]{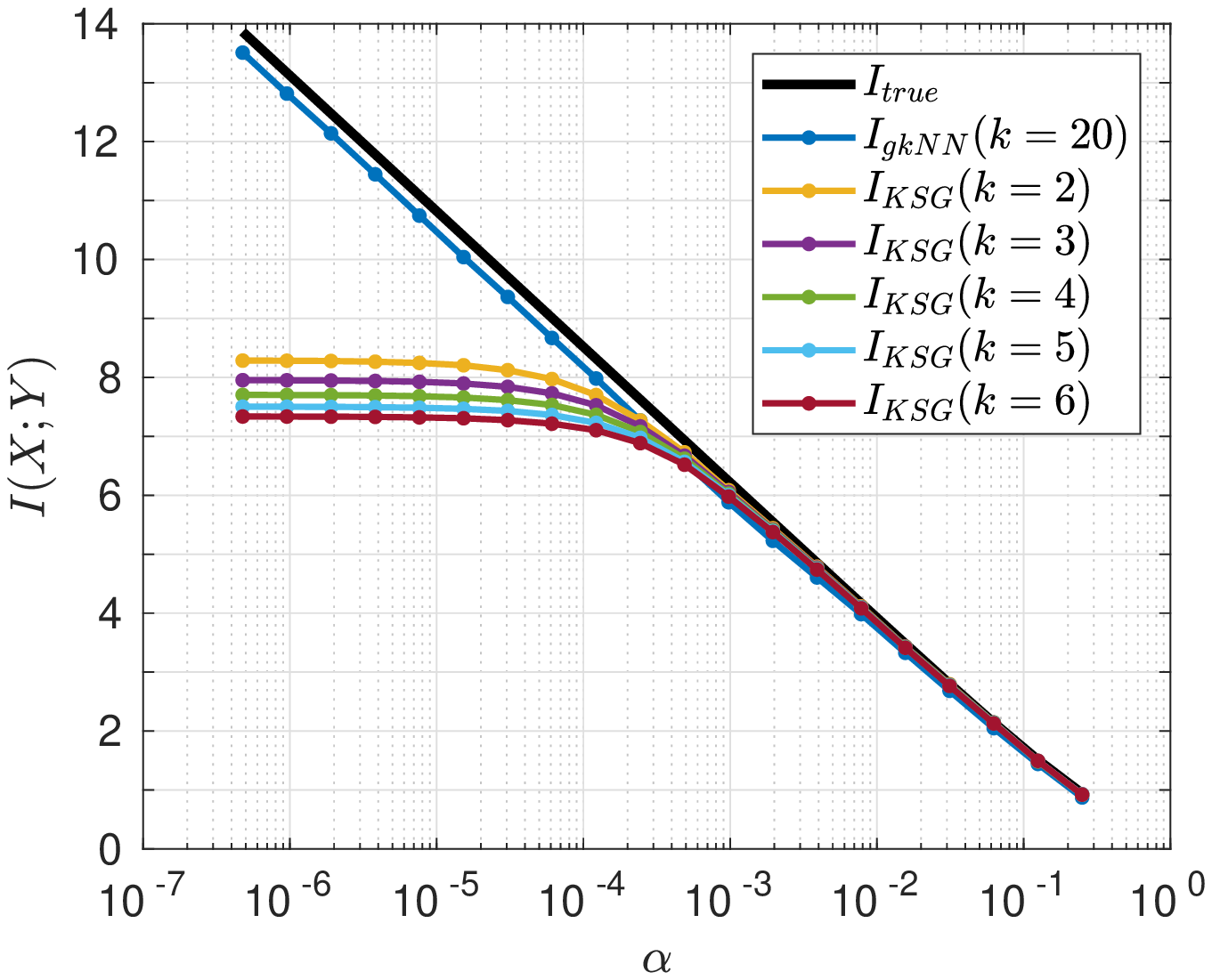}
    \caption{}
    \label{fig:alphavaryUnif2d}
  \end{subfigure}
  \begin{subfigure}[b]{.37\textwidth}
    \includegraphics[width=\textwidth]{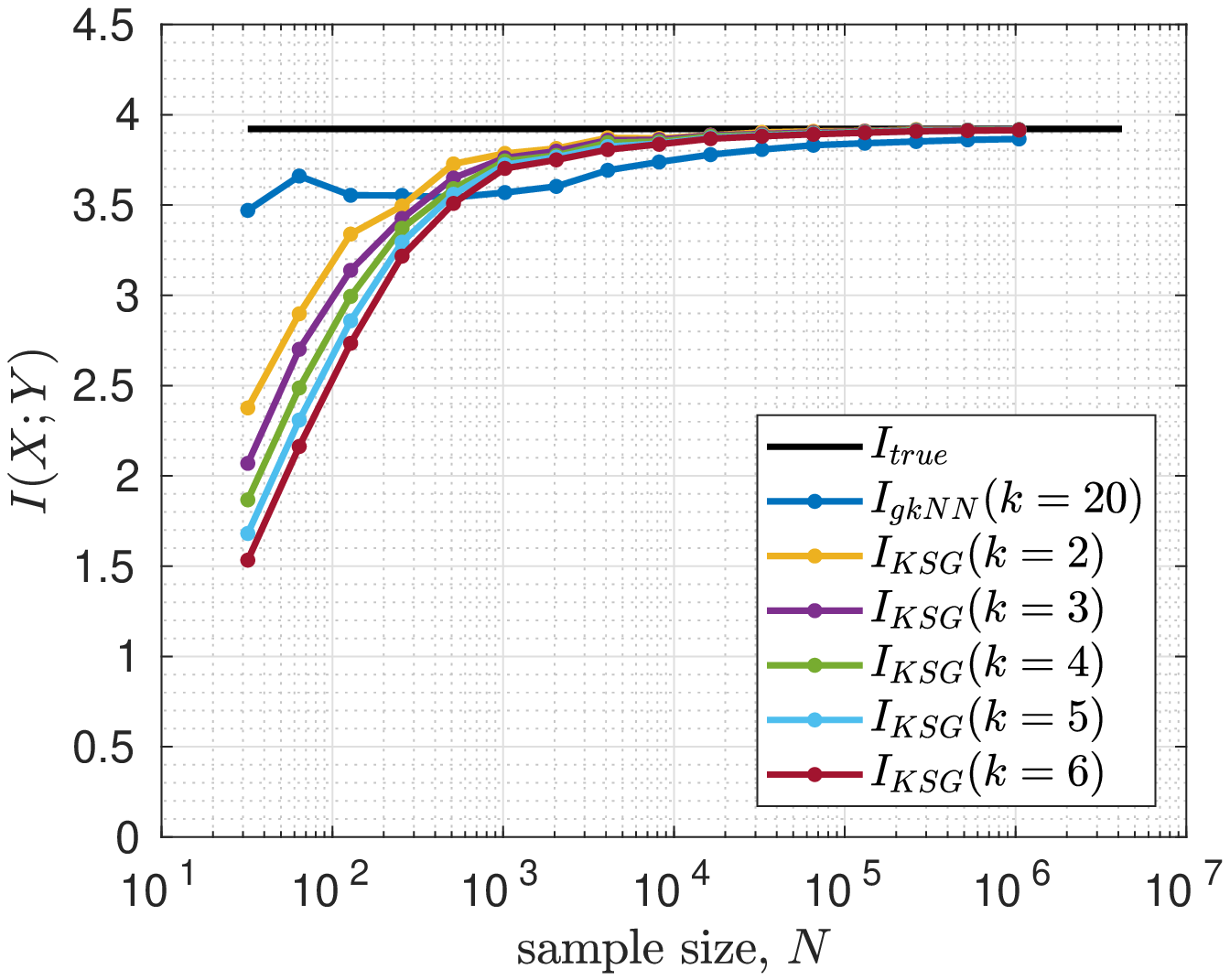}
    \caption{}
    \label{fig:NvaryUnif2d}
  \end{subfigure}\\
  \begin{subfigure}[b]{.37\textwidth}
    \includegraphics[width=\textwidth]{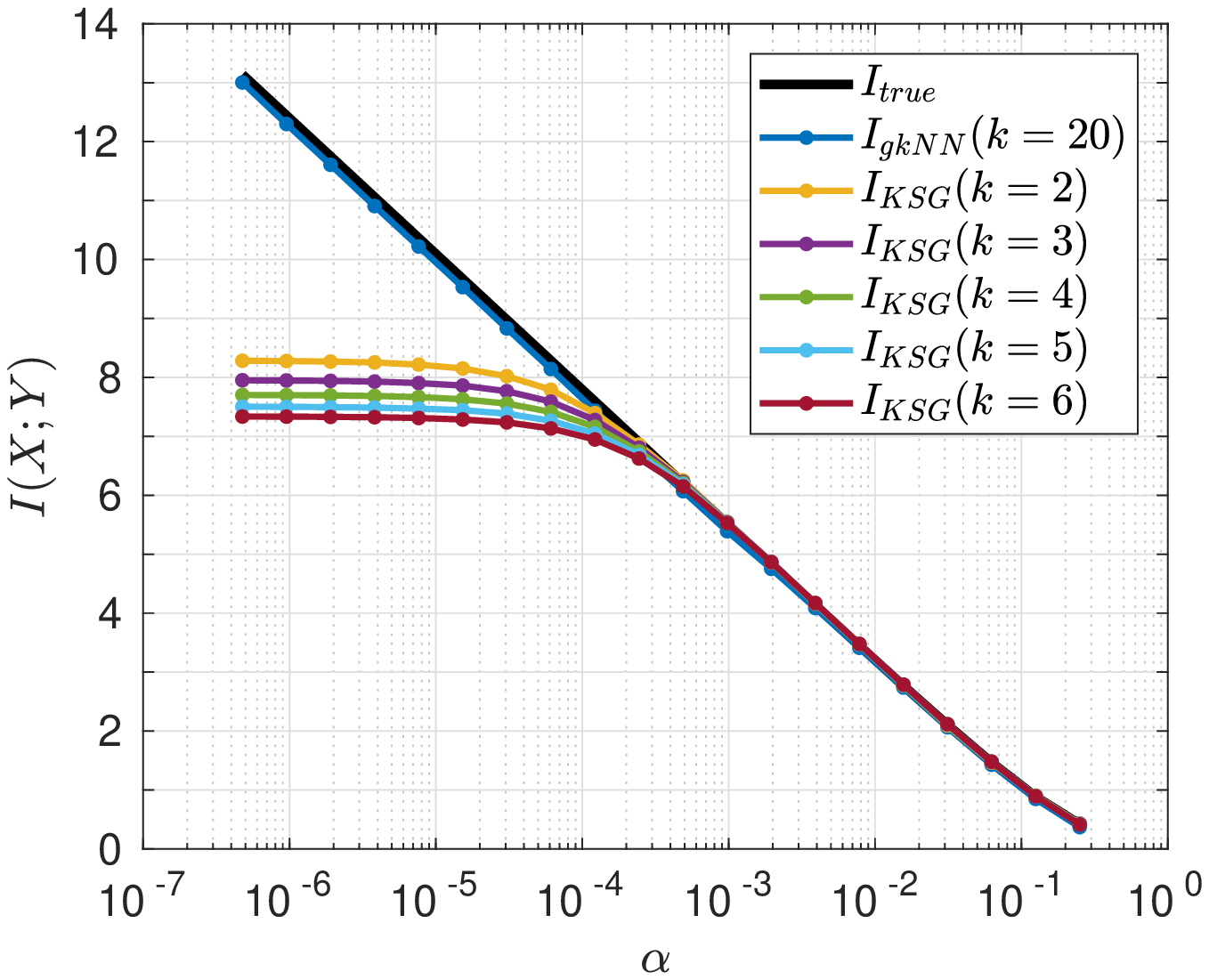}
    \caption{}
    \label{fig:alphavaryNormal2d}
  \end{subfigure}
  \begin{subfigure}[b]{.37\textwidth}
    \includegraphics[width=\textwidth]{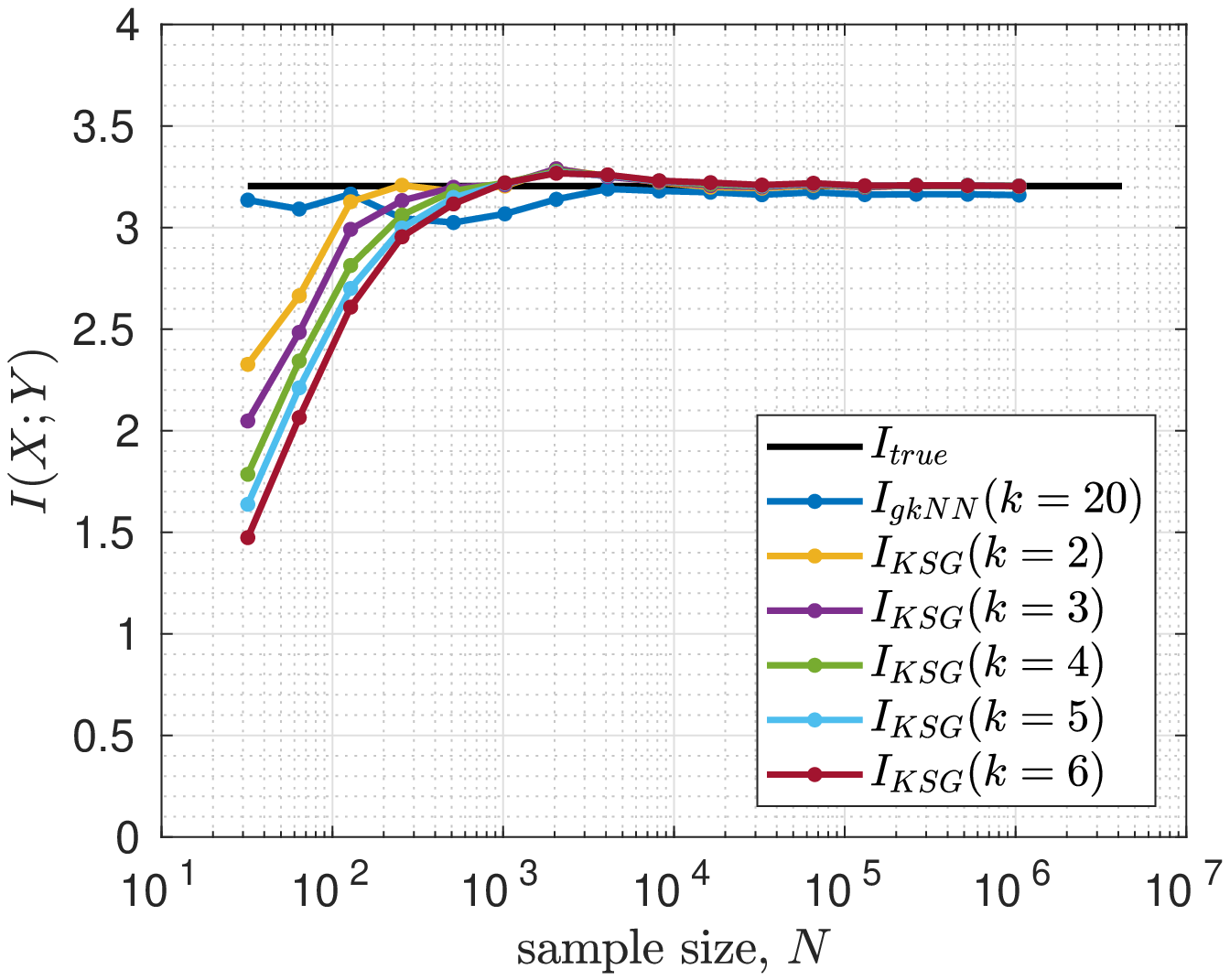}
    \caption{}
    \label{fig:NvaryNormal2d}
  \end{subfigure}\\
  \begin{subfigure}[b]{.37\textwidth}
   \includegraphics[width=\textwidth]{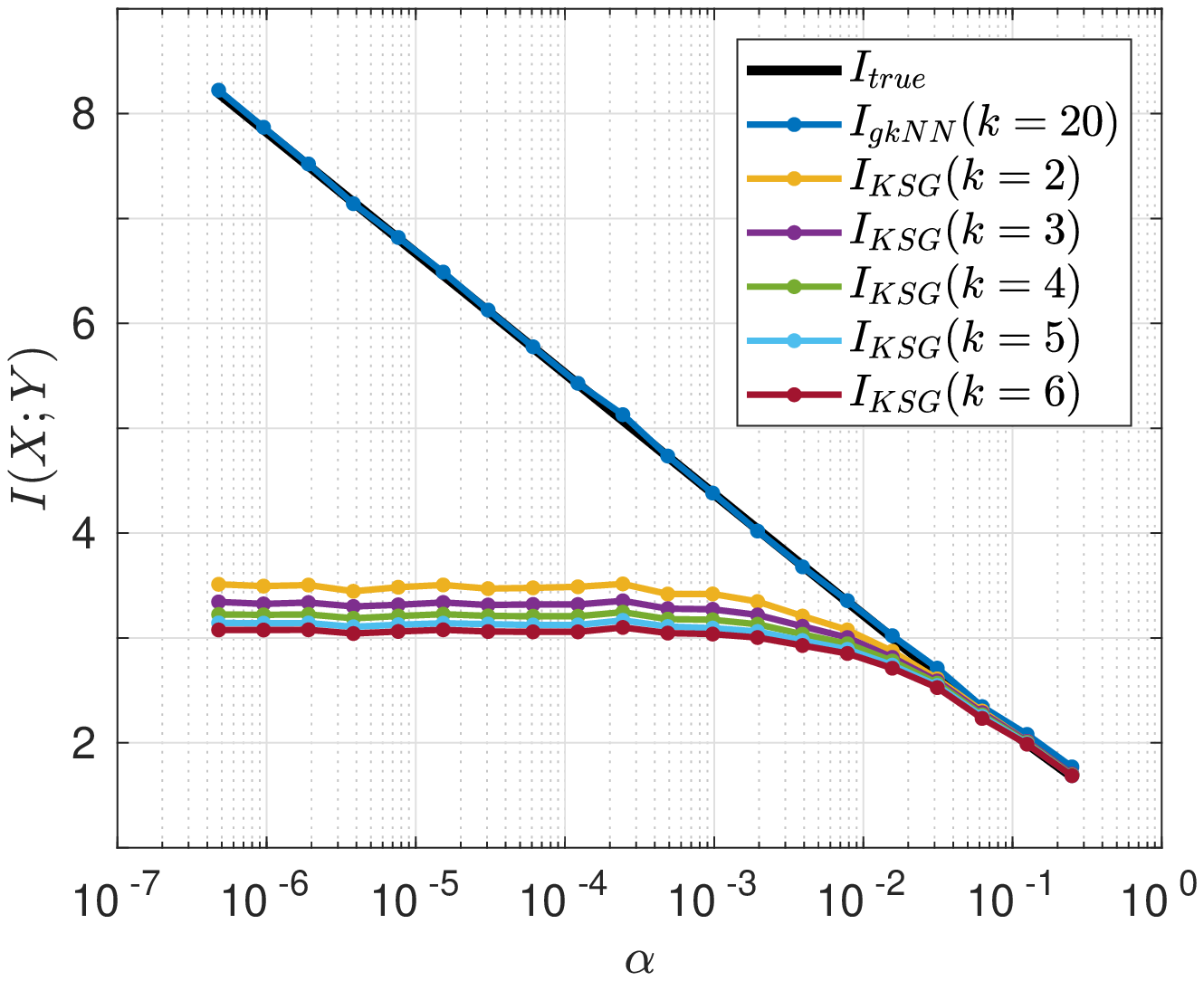}
    \caption{}
    \label{fig:alphavaryNormal4d}
  \end{subfigure}
    \begin{subfigure}[b]{.37\textwidth}
    \includegraphics[width=\textwidth]{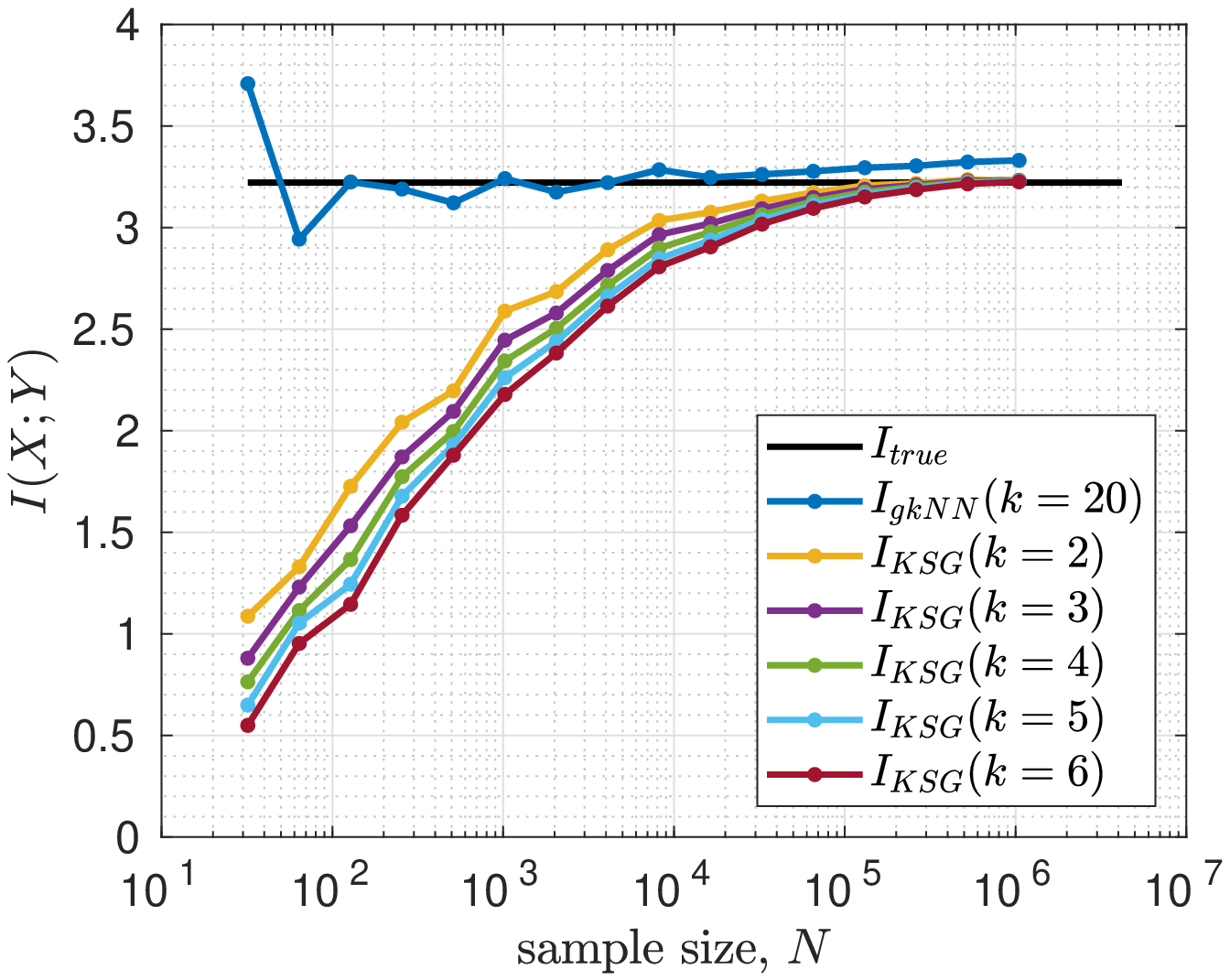}
    \caption{}
    \label{fig:NvaryNormal4d}
  \end{subfigure}
   \caption{A comparison of the KSG estimator with the g-k\textcolor{black}{nn} estimator on samples from the three families of variables. The top row of figures correspond to variables from family 1, the middle row to variables from family 2, and the bottom row to variables from family 3. In Figs. (a), (c), and (e), the sample size $N$ is fixed at $10^4$ and the thickness parameter $\alpha$ of each family is varied. On the right the thickness parameter of each family is fixed at $\alpha=.01$ and $N$ is allowed to vary. For each value of $\alpha$ (left) or $N$ (right) one sample of the joint random variable of size $N$ is drawn and both the KSG and g-k\textcolor{black}{nn} estimates are performed on the same sample.}
  \label{fig:results}
\end{figure*}

There are two \textcolor{black}{ways} in which the $k$-neighborhoods in these examples can get stretched. One way is that $\alpha$ gets small while $N$ stays fixed. The other is that the local neighborhoods determined by the $k$ nearest neighbors get larger. This occurs when $\alpha$ is small but fixed, and $N$ is decreased, because the sample points become more spread out.

Figure~\ref{fig:results} shows the results of KSG and svd estimates on samples of each of the three types of variables. In the figures on the left $N$ is fixed at $10^4$ and $\alpha$ varies. In the figures on the right $N$ is allowed to vary but $\alpha$ is fixed at $1/100$. The top row of figures correspond to samples from Family 1, the middle row to Family 2, and the bottom row to Family 3. For each value of $\alpha$ or $N$, one sample of size $N$ of the joint random variable is created and used by both estimators.

For each of the g-k\textcolor{black}{nn} estimates we have used $k=20$. The value $k=20$ was chosen because $k$ should be small enough to be considered a local estimate, but large enough that the svd of the $(k+1)\times d$ matrix of centered data should give good estimates of the directions and proportions of stretching. The value $k=20$ is chosen because it seems to balance these criteria, but no attempt has been made at optimizing $k$. Furthermore, in principle $k$ should depend on the dimension of the joint space because the number of axes of an ellipsoid is equal to the dimension of the space. Therefore a better estimate might be obtained by using a larger $k$ for Family 3 than the value of $k$ used in Families 1 and 2.

For the KSG estimator, $k$ is allowed to vary between 2 and 6. The value $k=1$ is excluded because it is not used in practice due to its large variance. The most common choices of $k$ are between 4 and 8~\cite{gao2017demystifying}. In each numerical example in this paper, however, the KSG estimates become progressively worse as $k$ increases, so that we omit the larger values of $k$ in order to better present the best estimates of KSG.

In each of Figs.~\ref{fig:alphavaryUnif2d},~\ref{fig:alphavaryNormal2d}, and~\ref{fig:alphavaryNormal4d} the true value of $I(X;Y)$ increases asymptotically like $\log(\alpha)$ as $\alpha\to0$. The KSG estimates do well for larger $\alpha$ but level off at a threshold that depends on the family\textcolor{black}{, indicating that local geometry is a likely cause of the bias}. The g-k\textcolor{black}{nn} estimator\textcolor{black}{, in contrast} keeps increasing like $\log(\alpha)$ as $\alpha$ gets smaller\textcolor{black}{, suggesting that the adaptations to traditional knn methods allow the g-knn method to adapt to the changing local geometry}.

In each of Figs.~\ref{fig:NvaryUnif2d},~\ref{fig:NvaryNormal2d}, and~\ref{fig:NvaryNormal4d}, the true value is constant for all $N$ since it is determined by the distribution, which depends on $\alpha$ alone. For large $N$\textcolor{black}{,} KSG outperforms g-k\textcolor{black}{nn} since it is asymptotically unbiased \cite{gao2017demystifying}. As sample size is decreased, however, at some point the KSG estimate becomes progressively more biased. This is because fewer samples means that the data points are more spread apart, which, with a fixed $k$, means the $k$-neighborhoods are larger. Therefore for some $N$ the neighborhoods effectively span the width of the distribution in the thinnest direction. The g-k\textcolor{black}{nn} estimate on average stays near the true value as $N$ decreases, although its variance seems to increase. It stays close to the true value because it is able to adapt to the stretched $k$-neighborhoods by using thinner ellipsoids.

\subsection{A more complex example}\label{sec:henon}
\textcolor{black}{Note that Figure~\ref{fig:results} should not be interpreted as a direct comparison of KSG and g-knn because clever preprocessing could be applied to data sampled from each family to remove the flatness from the local geometry.} A more complex example\textcolor{black}{, which would likely resist efforts to preprocess the data to counteract the effects of local geometry,} is provided by a 4d system consisting of coupled H\'{e}non maps. In this system both $X$ and $Y$ are \textcolor{black}{2d} and $Y$ is coupled to $X$, but not vice-versa, so that $X$ can be thought of as driving $Y$. The purely deterministic system approaches a measure 0 attractor so that $I(X;Y)$ would not be defined without the addition of noise, which is added to the $X$ variable, and reaches the $Y$ variable through the coupling.

\begin{description}
  \item{\textcolor{black}{\textbf{Family 4: }}}The system is
\begin{align}\label{eq:henondefbegin}
  X_{1,n+1} &= a-X_{1,n}^2 + bX_{2,n} + \eta_1\\
  X_{2,n+1} &= X_{1,n} +\eta_2\\
  Y_{1,n+1} &= a - (c X_{1,n}Y_{1,n} + (1-c) Y_{1,n}^2)+bY_{2,n}\\
  Y_{2,n+1} &= Y_{1,n} \\
  \eta_1&\sim  Unif(-\alpha,\alpha)\\
  \eta_2&\sim  Unif(-\alpha,\alpha),\label{eq:henondefend}
\end{align}
where $a=1.2$, $b=0.3$, and the coupling coefficient is $c=0.8$. When $\alpha=0$ the system is the same coupled H\'{e}non map described in a number of studies~\cite{kreuz2007measuring} except that $a$ is reduced from the usual 1.4 to 1.2 so that noise can be added without causing trajectories to leave the basin of attraction. In these studies it is noted that a coupling coefficient of $c=0.8$ results in identical synchronization so that in the limit of large $n$, $X_{1,n}=Y_{1,n}$ and $X_{2,n}=Y_{2,n}$, implying that the limit set is contained in a 2d manifold. Thus, in the long term, and as $\alpha\to0$, samples from this stochastic process should lie near a 2d submanifold of $\R^4$.
\end{description}

\begin{figure}
  \begin{subfigure}[b]{.9\columnwidth}
    \includegraphics[width=\linewidth]{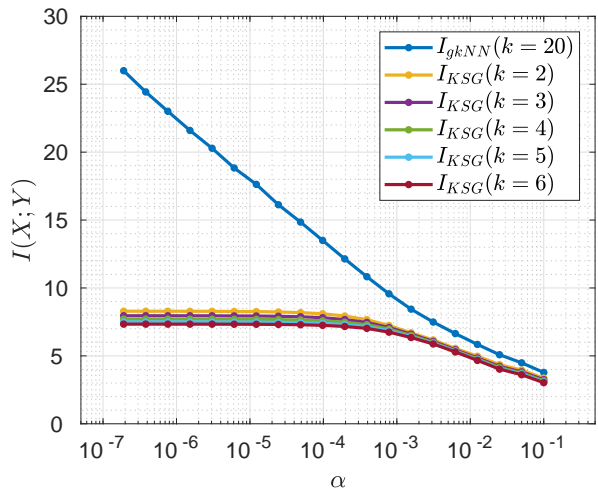}
    \caption{}
    \label{fig:alphavaryHenon}
  \end{subfigure}

  \vspace{1em}
  \begin{subfigure}[b]{.9\columnwidth}
    \includegraphics[width=\linewidth]{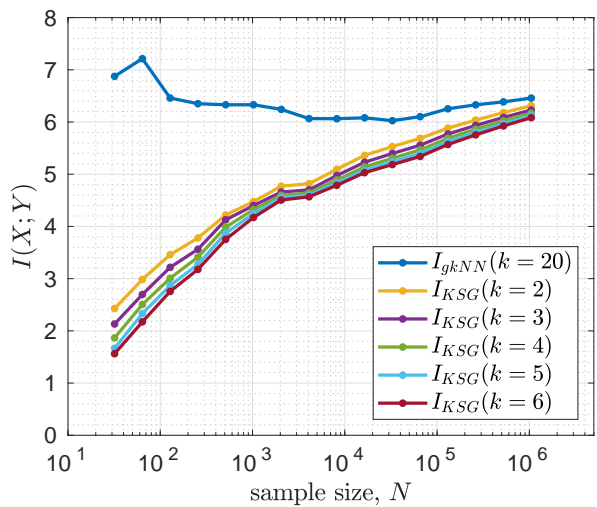}
    \caption{}
    \label{fig:NvaryHenon}
  \end{subfigure}
  \caption{A comparison of the KSG estimator with the g-k\textcolor{black}{nn} estimator for the stochastic coupled H\'{e}non map described in Eqs.\eqref{eq:henondefbegin} through~\eqref{eq:henondefend}. In (a) the sample size is fixed at $10^4$ while $\alpha$ is allowed to vary. In (b), $\alpha$ is fixed at $\alpha=.01$ and sample size varies.}
  \label{fig:Henon}
\end{figure}

\textcolor{black}{It is important to note that the noise is introduced dynamically, and transformed by a nonlinear transformation on each time step. The samples lie near the H\'{e}non attractor embedded in the 2d submanifold. If one thinks of the data as H\'{e}non attractor plus noise, then the noise depends on $X$ and $Y$ in a nonlinear manner, and produces heterogeneous local geometries.}

Figure~\ref{fig:Henon} compares the estimates given by the g-k\textcolor{black}{nn} estimator and the KSG estimator. The exact value of $I(X;Y)$ is unknown, but there are qualitative differences between the performance of the two estimators. Fig.~\ref{fig:alphavaryHenon} compares the estimates as $\alpha$ is decreased, in which case we expect that the true value of $I(X;Y)$ increases unboundedly, a behavior captured by the g-k\textcolor{black}{nn} estimator, but not the KSG estimator, which appears asymptotically constant as $\alpha\to 0$.

In Fig.~\ref{fig:NvaryHenon}, the sample size is varied, which has no effect on the underlying mutual information, $I(X;Y)$. Compared to the KSG estimates, the g-k\textcolor{black}{nn} estimates are relatively constant as $N$ is varied. Of particular practical importance is the behavior as $N$ is decreased, which seems to introduce a lot of negative bias into the KSG estimate. The g-k\textcolor{black}{nn} estimates increase slightly as $N$ decreases, which might indicate a slight positive bias, but it is difficult to tell with only one sample per \textcolor{black}{plotted} point. This issue might be more thoroughly investigated using the mean and variance of a \textcolor{black}{large} number of estimates for each value of $N$.

\section{Discussion}

A common strategy in $k$nn estimation of differential entropy and mutual information is to use local data to fit volume elements. The use of geometrically regular volume elements \textcolor{black}{requires minimal local data,} so that the volume elements remain as local\textcolor{black}{ized} as possible. This paper introduces the notion of a  g-k\textcolor{black}{nn} estimator, which uses slightly more data points to fit local volume elements in order to better model the local geometry of the underlying measure.

As an application, this paper derives a g-k\textcolor{black}{nn} estimator of mutual information, inspired by a consideration of the local geometry of dynamical systems attractors. A common feature of dissipative systems and systems with competing time scales is that their limit sets lie in a lower dimensional attractor or manifold. Locally the geometry is typically characterized by directions of maximal stretching and compression, which are described quantitatively by the Lyapunov spectrum. Ellipsoids are used for local volume elements because they capture the directions of stretching and compression without requiring large amounts of local data to fit.

It might be noted that the ellipsoids are simply spheres in the Mahalanobis distance determined by the local data \cite{mahalanobis1936generalised}. The metric that is used to define the sphere\textcolor{black}{s} in the g-k\textcolor{black}{nn} estimate, however, varies from neighborhood to neighborhood \textcolor{black}{. This behavior} is very different from many other \textcolor{black}{knn} estimators of pdfs \textcolor{black}{where the spheres} are determined by a global metric (typically defined by a $p$-norm). In this perspective, g-k\textcolor{black}{nn} methods use data to learn both local metrics and volumes, and hence a local geometry, justifying the use of the name g-k\textcolor{black}{nn}.

\textcolor{black}{The numerical examples suggest that when it is not possible to preprocess the data to add thickness to the local geometry,} the g-k\textcolor{black}{nn} estimator outperforms KSG as the underlying measure becomes more thinly supported. \textcolor{black}{T}he g-k\textcolor{black}{nn} estimator also outperforms KSG as sample size decreases, a result that is particularly promising for applications in which the number of data points is limited.
However, unlike the \textcolor{black}{Kozachenko-Leonenko} estimator of differential entropy and the KSG estimator of mutual information, the g-knn estimator as based on ellipses developed in this paper is not asymptotically unbiased. In our future work we hope to remove this asymptotic bias in a manner analogous to Ref.~\cite{singh2003nearest} to gain greater accuracy for both low and high values of $N$.

There are also other descriptions of local geometry that suggest alternative g-knn methods. For instance, there are nonlinear partition algorithms such as OPTICS~\cite{ankerst1999optics}, which are based on data clustering. Also, entropy is closely related to recurrence, which is suggestive of alternative g-knn methods based on the detection of recurrence structures~\cite{beim2016optimal}.

\begin{acknowledgements}
This work was funded by ARO Grant No. W911NF-12-1-0276, Grant No. N68164-EG, and ONR Grant No. N00014-15-1-2093.
\end{acknowledgements}


%

\end{document}